\documentclass[11pt]{amsart}
\usepackage{enumerate}
\usepackage{geometry}                
\usepackage{graphicx}
\usepackage{amssymb}
\usepackage{epstopdf}
\usepackage{geometry} 
\usepackage{graphicx}
\usepackage{amsmath} 
\usepackage{stackengine}
\usepackage{romanbar}
\numberwithin{equation}{section}
\usepackage{commath}
\usepackage{cleveref}
\usepackage[normalem]{ulem}
\usepackage{bbm}
\usepackage{subfig}
\usepackage{float}
\usepackage{url}
\usepackage{mathtools}
\usepackage{amsthm} 
\usepackage[disable]{todonotes} 
\theoremstyle{plain}
\newtheorem{theorem}{Theorem}[section]
\newtheorem*{theorem*}{Theorem}
\newtheorem{lemma}[theorem]{Lemma}

\newtheorem{proposition}[theorem]{Proposition}

\theoremstyle{remark}
\newtheorem{remark}{Remark}
\theoremstyle{remark}

\theoremstyle{definition}
\newtheorem{definition}[theorem]{Definition}
\allowdisplaybreaks
\title[$BV$ blowup for $p$-system]{Finite time BV blowup for Liu-admissible solutions to $p$-system via computer-assisted proof}
\author[Sam G. Krupa]{Sam G.  Krupa}
\address[Sam G. Krupa]{\newline Max Planck Institute for Mathematics in the Sciences, \newline 04103 Leipzig, Germany}
\email{Sam.Krupa@mis.mpg.de}

\thanks{Funded by the German Research Foundation (DFG) project number 525859002. The author would also like to thank L\'{a}szl\'{o} Sz\'{e}kelyhidi, Jr. for helpful discussions.}

\date{\today}                                           

\begin{document}
\keywords{Systems of conservation laws, one space dimension, genuine nonlinearity, isentropic gas dynamics, Liu entropy condition, entropy conditions, blowup, bounded variation, non-uniqueness, convex integration, differential inclusion.}
\subjclass[2020]{Primary 35L65; Secondary 35B44, 35D30, 35A02, 35L45, 76N15, 49K21}
\begin{abstract}
   In this paper, we consider finite time blowup of the $BV$-norm for exact solutions to genuinely nonlinear hyperbolic systems in one space dimension, in particular the $p$-system. We consider solutions verifying shock admissibility criteria such as the Lax \emph{E}-condition and the Liu \emph{E}-condition. In particular, we present Riemann initial data which admits infinitely many bounded solutions, each of which experience, not just finite time, but in fact \emph{instantaneous} blowup of the $BV$ norm. The Riemann initial data is allowed to come from an open set in state space. Our method provably does not admit a strictly convex entropy. 
   
   The main results in this article compare to Jenssen [{\em SIAM J. Math. Anal.}, 31(4):894--908, 2000], who shows $BV$ blowup for bounded solutions, or alternatively, blowup in $L^\infty$, for an artificial $3\times 3$ system which is not genuinely nonlinear. Baiti-Jenssen [{\em Discrete Contin. Dynam. Systems}, 7(4):837--853, 2001] improves upon this Jenssen result and can consider a genuinely nonlinear system, but then the blowup is only in $L^\infty$ and they cannot construct bounded solutions which blowup in $BV$. Moreover, their system is non-physical and provably does not admit a global, strictly convex entropy. Our result also shows sharpness of the recent Bressan-De Lellis result [{\em Arch. Ration. Mech. Anal.}, 247(6):Paper No. 106, 12, 2023] concerning  well-posedness via the Liu \emph{E}-condition. The proof of our theorem is computer-assisted, following the framework of Sz\'{e}kelyhidi [{\em Arch. Ration. Mech. Anal.}, 172(1):133--152, 2004]. Our code is available on the GitHub.
\end{abstract}

\maketitle
\tableofcontents

\section{Introduction}

Consider the $p$-system, which is a $2\times 2$ system of conservation laws in one space dimension,
\begin{align}\label{psystem}
    \begin{cases}
    \partial_t v+\partial_x u=0,\\
    \partial_t u +\partial_x p(v)=0,
    \end{cases}
\end{align}
for the unknown functions $v,u\colon \mathbb{R}\times[0,\infty)\to\mathbb{R}$ and where $t>0$ and $x\in\mathbb{R}$. We consider functions $p\colon\mathbb{R}\to\mathbb{R}$ such that $p'>0$ and $p''<0$. The $p$-system is a model of isentropic gas dynamics in Lagrangian coordinates. Remark that in this setting, the $p$-system is strictly hyperbolic and genuinely nonlinear. 

We will consider the Cauchy problem, with initial data
\begin{align}\label{init}
 \begin{cases}
    v(x,0)=v^0(x),\\
    u(x,0)=u^0(x),
    \end{cases}
\end{align}
for functions $v^0,u^0\colon\mathbb{R}\to\mathbb{R}$. 

In the introduction to \cite{MR3914001}, Bressan-Chen-Zhang write ``For hyperbolic systems of conservation laws in one space dimension \cite{MR1816648,MR1301779,MR1723032,MR1912206}, a major remaining open problem is whether, for large $BV$ initial data, the total variation of entropy-weak solutions remains uniformly bounded or can blow up in finite time.''

This paper is dedicated to the proof of the following theorem.

\begin{theorem}[Main theorem -- BV blowup for the $p$-system]\label{main_theorem}
There exists a smooth, strictly increasing, and strictly concave function $p\colon\mathbb{R}\to\mathbb{R}$ such that for the $p$-system \eqref{psystem} with this $p$, the following holds:
 There exists open sets $\mathcal{U}_L,\mathcal{U}_R\subset\mathbb{R}$ such that for all Riemann initial data
 \begin{align}\label{riemann_initial}
(v^0(x),u^0(x))=
\begin{cases}
(v_L,u_L) & \mbox{ for } x<0\\
(v_R,u_R) & \mbox{ for } x>0,
\end{cases}
\end{align}
with $(v_L,u_L)\in \mathcal{U}_L$ and $(v_R,u_R)\in \mathcal{U}_R$, there exists infinitely many different solutions $(v,u)\in L^\infty (\mathbb{R}\times[0,\infty))$ to the $p$-system \eqref{psystem} with the property that for every time $t>0$, the function $(v(\cdot,t),u(\cdot,t))$ has infinite total variation. Moreover, each solution will only take 5 different values, no two of which are connected by a shock. Thus, each solution vacuously verifies the Lax \emph{E}-condition and Liu \emph{E}-condition. 
\end{theorem}

\Cref{main_theorem} is the first result we know of which gives finite time $BV$ blowup of bounded solutions for a genuinely nonlinear system with admissible shocks.

Broadly speaking, the technique of proof is convex integration, which is a method of constructing highly oscillatory solutions to partial differential equations (PDEs).

In the setting of hyperbolic systems of conservation laws (such as \eqref{psystem}), weak (in the sense of distribution/measure) solutions seem to be the natural class of solutions to consider. In fact, even in the case of scalar equations in one space dimension, simple examples show that $C^\infty$ initial data may evolve discontinuities (``shocks'') in finite time. Such weak solutions are easily seen to be highly non-unique. The hope then is that some mix of various selection criteria will select for a unique admissible (``physical'') solution for each fixed initial data. For example, to determine if a particular shock discontinuity is admissible, selection criteria such as the Lax \emph{E}-condition or the Liu \emph{E}-condition are employed (see \Cref{sec:shock_conds}, below). Another direction is to consider a strictly convex entropy functional and an associated entropy inequality applied to the weak solution itself (see \cite{dafermos_big_book}, as well as the theory of shifts of Vasseur \cite{MR2508169}, recently culminating in \cite{MR4487515}).

To be more precise, we define a weak solution to \eqref{psystem}, \eqref{init} as a 2-tuple $(v,u)$ which verifies 
\begin{equation}
\int_{0}^{\infty} \int_{\mathbb{R}} \left[ \partial_t \Phi \cdot U + \partial_x \Phi \cdot f(U) \right] \, dxdt + \int_{\mathbb{R}} \Phi(x, 0)\cdot U^0(x) \, dx = 0,
\end{equation}
for every Lipschitz test function \( \Phi: \mathbb{R} \times [0,\infty) \rightarrow \mathbb{R}^{2} \) with compact support, and where we write $U\coloneqq (v,u)$, $U^0\coloneqq (v^0,u^0)$  and $f(v,u)\coloneqq (u,p(v))$.

In the context of $2\times 2$ hyperbolic systems, \emph{strictly hyperbolic} means that the eigenvalues \( \lambda_1(U) <\lambda_2(U) \) of the Jacobian \( Df(U) \) are real and distinct for all  \( U\in\mathbb{R}^2 \). Let \( r_1(U), r_2(U) \) denote the corresponding eigenvectors, normalized with unit length. The \( k \)-th characteristic field is said to be \emph{genuinely nonlinear} if
\[
\nabla \lambda_k(U) \cdot r_k(U) \neq 0 \text{ for all } U,
\]
while it is \emph{linearly degenerate} if
\[
\nabla \lambda_k(U) \cdot r_k(U) = 0.
\]

If both characteristic fields are genuinely nonlinear we say the system is a \emph{genuinely nonlinear system}.

\subsection{Finite time blowup}

The finite-time $BV$ blowup of bounded, entropy-admissible solutions to a \emph{genuinely nonlinear} hyperbolic system, in one space dimension, is an important and long-standing open problem in the field. Let us recall some of the celebrated works in this direction.

\subsubsection{Joly-M\'{e}tivier-Rauch (1994) and Young (1999): magnification}
The works \cite{MR1677943,MR1272766} give examples of $3\times3$ systems which show arbitrary, unbounded \emph{magnification} of the total variation or $L^\infty$-norm of solutions. We remark that the work  \cite{MR1272766} considers genuinely nonlinear systems, while the work \cite{MR1677943} considers only systems where all characteristic fields are linearly degenerate.

\subsubsection{Jenssen (2000) and Baiti-Jenssen (2001): blowup for artificial systems}
In \cite{MR1752421}, Jenssen gives examples of $3\times 3$ systems of strictly hyperbolic conservation laws which exhibit solutions for which the $L^\infty$-norm of the solution blows up in finite time, or the total variation of the solution explodes in finite time while the $L^\infty$-norm of the solution remains bounded. These are the first examples we are aware of where the $L^\infty$-norm or $BV$-norm truly becomes \emph{infinite} in \emph{finite time} (in contrast with \cite{MR1677943,MR1272766}, discussed above). However, the systems constructed by Jenssen are linearized in the sense that the first and third characteristic fields are linearly degenerate.

In \cite{MR1849664}, Baiti-Jenssen construct an artificial $3\times 3$ system, not motivated by physics, which is genuinely nonlinear in all characteristic fields, and they exhibit piecewise-constant initial data, with three shocks, and which gives finite time blowup in $L^\infty$. Baiti-Jenssen works on perturbations of systems introduced in the earlier work \cite{MR1752421} by Jenssen, which considers $3\times 3$ systems where the first and third characteristic fields are linearly degenerate.   The Baiti-Jenssen construction is based upon carefully building a particular wave pattern with a ``ping-pong'' behavior where infinitely many shock waves are produced in finite time (and Riemann problems are solved at each interaction). But they are not able to construct solutions which blow up in $BV$ while remaining bounded. In their situation, the geometry is too delicate to execute such a construction. They also prove that their $3\times 3$ artificial system does not admit any global, strictly convex entropy. However, they utilize wave families verifying the Lax \emph{E}-condition and Liu \emph{E}-condition \footnote{Baiti-Jenssen do not explicitly show that their wave families verify the Lax \emph{E}-condition and Liu \emph{E}-condition. We give a short proof of these facts in the appendix (\Cref{sec:laxliu_appendix}).}.

In contrast, our result (\Cref{main_theorem}) gives exact solutions with simple Riemann initial data which exhibit \emph{instantaneous} blowup in $BV$. Unlike Baiti-Jenssen \cite{MR1849664}, we are able to work with a physical system (the $p$-system), and we require only two wave families as opposed to three. Our solutions contain no shocks, and thus verify all shock-based admissibility criteria including the Lax \emph{E}-condition (see \eqref{lax_cond}) and Liu \emph{E}-condition (see \eqref{liu_cond}). Similarly to Baiti-Jenssen, our techniques  provably do not admit a convex entropy.  This is due to our use of $T_5$ configurations (explained below in \Cref{sec:sol_by_convex_int}). See Johansson-Tione \cite{https://doi.org/10.48550/arxiv.2208.10979} for a proof that a $T_5$ does not exist for the $p$-system if a strictly convex entropy is assumed. With the assumption of a strictly convex entropy, see also Lorent-Peng \cite{MR4144350} for a proof that the $p$-system does not admit $T_4$ configurations. In a recent and closely related work, the author and Sz\'{e}kelyhidi prove the nonexistence of $T_4$ configurations for a large class of hyperbolic systems with a strictly convex entropy \cite{2022arXiv221114239K}.

\subsubsection{Bressan-Chen-Zhang (2018): perturbation of wave speeds for the $p$-system}
In \cite{MR3914001}, the authors study possible mechanisms for $BV$ blowup in the $p$-system, in the genuinely nonlinear regime $p''\neq 0$. They consider piecewise-smooth approximate solutions. If wave speeds are slightly perturbed, then this changes the order in which waves interact. Thus, from the same initial data, approximate solutions are constructed with either bounded total variation, or total variation which blows up in finite time. The key is the interaction pattern of the waves. The authors write \cite[p.~1243]{MR3914001} ``Although our solutions are not exact, because some errors occur in the wave speeds, they possess all the qualitative properties known for exact solutions. The present analysis thus provides some indication that finite time blowup of the total variation might be possible, for the $p$-system.'' For the associated exact solutions, their properties are unknown. 

\subsubsection{Bressan-De Lellis (2023): well-posedness via Liu \emph{E}-condition} For $n\times n$ hyperbolic systems of conservation laws, there is a well-known $L^1$-type semigroup of solutions which are assumed to have small $BV$ at each fixed time $t$. The semigroup coincides with the limits of front tracking approximations \cite{MR1816648}, limits of vanishing viscosity approximations \cite{MR2150387}, as well as solutions from the Glimm scheme \cite{MR0194770}. 

In \cite{MR4661213}, Bressan and De Lellis show that every weak solution taking values in the domain of this semigroup, and with shocks verifying the Liu \emph{E}-condition \eqref{liu_cond}, is actually equal to a semigroup trajectory and is thus unique. They make no assumption about genuine nonlinearity, and they \emph{do not assume the existence of a convex entropy.} Their only stipulation is that all points of approximate jump discontinuity verify the Liu \emph{E}-condition \eqref{liu_cond}.

\Cref{main_theorem} shows that this recent Bressan-De Lellis result is sharp. Our solutions vacuously verify the Liu \emph{E}-condition, but we are outside the $BV$ regime and our solutions are non-unique.

\subsubsection{The $L^2$ theory (2022)}

Recent work of the author, Chen, and Vasseur \cite{MR4487515} also considers the classical $L^1$ semigroup of solutions with small $BV$ (e.g. \cite{MR1816648,MR2150387,MR0194770}), and shows that these solutions are stable, and do not blowup, even when considered in the class of large $L^2$-perturbations. The work \cite{MR4487515} considers systems endowed with a strictly convex entropy, and assumes a trace condition on the solutions. However, it is interesting to study the stability of $BV$ data in an even broader class. In the present paper, we expand the class of data from \cite{MR4487515} to the wider class of $L^\infty$ solutions verifying the Lax \emph{E}-condition and Liu \emph{E}-condition.

\subsubsection{Chiodaroli-De Lellis-Kreml (2015): convex integration solutions for 2-D Riemann problems} In \cite{multi_d_illposed}, Chiodaroli-De Lellis-Kreml consider the isentropic compressible Euler system in \emph{two space dimensions}. They show that bounded solutions to a one-dimensional Riemann problem can be constructed via convex integration (for related work, see also \cite{MR4385531}). Their result gives instantaneous blowup of the local $BV$-norm. It is interesting to note that the solutions constructed via convex integration to the 1-D Riemann problems are truly two-dimensional in nature. The solutions constructed by Chiodaroli-De Lellis-Kreml verify an entropy inequality for a strictly convex entropy. The work utilizes the framework for convex integration developed for the \emph{incompressible theory} (see \cite{MR2564474,MR2600877}) and thus does not extend to the systems of conservation laws in one space dimension, where there is no corresponding incompressible theory. 

\subsection{Genericity and a question of Dafermos}
In the recent high-level survey article of Dafermos, he questions whether convex integration solutions and finite time blowup solutions (such as from \Cref{main_theorem}) are generic phenomena \cite[p.~484-485]{MR4627977}.

Concerning convex integration solutions for compressible equations, Dafermos writes  \cite[p.~-485]{MR4627977} ``A possible explanation is that we are missing the proper admissibility condition on weak solutions. Alternatively, it may turn out that exotic solutions are not observed in nature because they are not generic.'' 

In our context, we are able to give an answer to Dafermos's question: indeed, we are able to construct solutions with Riemann initial data coming from an open set in the state space, and which exhibit highly irregular behavior.

\subsection{Method of proof: convex integration}

We use the Kirchheim-M\"{u}ller-\v{S}ver\'{a}k framework \cite{MR2008346} to convert the question of solutions to the PDE \eqref{psystem} to the question of solutions to a differential inclusion of the form $D\psi\in\mathcal{K}$ for a certain constitutive set $\mathcal{K}\subset \mathbb{R}^{2\times 2}$ (this will be explained in more detail in \Cref{sec:sol_by_convex_int}, below). Then the question of solutions to \eqref{psystem} becomes a question of the rank-one convex geometry of the set $\mathcal{K}$, in particular the rank-one convex hull $\mathcal{K}^{rc}$.

A first step in the study of the rank-one convexity of $\mathcal{K}$ is to ask if $\mathcal{K}$ has rank-one connections, i.e. two matrices $A,B\in\mathcal{K}$ such that $\mathrm{rank}(A-B)=1$. If the set $\mathcal{K}$ has rank-one connections, then it is well-known that simple plane-wave solutions to the differential inclusion will exist. In our context, rank-one connections in $\mathcal{K}$ will correspond to shock solutions of the hyperbolic system \eqref{psystem} (see \Cref{sec:rank_one_shocks}, below).

\subsection{Large $T_5$ configurations}

A more interesting question is to study non-trivial (i.e., non-affine) solutions to the differential inclusion when $\mathcal{K}$ contains no rank-one connections. In particular, it is known that such solutions can exist even when the set $\mathcal{K}$ is small. In the work of Kirchheim-Preiss \cite[p.~100]{hab_thesis}, they give an example of a set $\mathcal{K}$ with 5 distinct points such that the corresponding differential inclusion admits non-trivial solutions. The work Kirchheim-Preiss \cite[p.~100]{hab_thesis} has been generalized by F\"{o}rster-Sz\'{e}kelyhidi \cite{MR3740399}, who create a general framework for the study of so-called ``non-rigid'' sets with 5 points and no rank-one connections. They study 5 point sets which are in so-called $T_5$ configuration (see \cite{MR1320538,MR2624766,MR1983780,MR2048569,MR2008346,MR2118899}, and \Cref{sec:TN_rank_one}, below) and develop a key generalization of this concept -- the \emph{large $T_5$ configuration}.  A crucial advantage to their characterization of non-rigid 5-element sets is that it provides for an algebraic criterion (see \Cref{alg_crit_theorem}, below) which can be implemented numerically without much effort, as we do in this paper, and which obviates the need for computing directly the rank-one convex hull.

With the correspondence between solutions to the differential inclusion and solutions to the system \eqref{psystem}, the F\"{o}rster-Sz\'{e}kelyhidi framework will give solutions to \eqref{psystem} which take exactly 5 values, no two of which can be connected by a shock. From this construction, our solutions will vacuously verify the Lax \emph{E}-condition (see \eqref{lax_cond}) and Liu \emph{E}-condition (see \eqref{liu_cond}).

The following exposition is based on \cite{MR3740399}: More precisely, F\"{o}rster-Sz\'{e}kelyhidi use the approach of \cite{MR1728376,MR1983780} and rely on \emph{in-approximations} (see \cite[p.~218]{MR864505}) to carry out the convex integration scheme.

We define the set  \( \Sigma \subset \mathbb{R}^{2 \times 2} \) as
\[
\Sigma \coloneqq \{X \in \mathbb{R}^{2 \times 2} : X \text{ is symmetric}\}.
\]

 We can now give the following definition and theorem:

\begin{definition}[In-approximations]\label{in_approx_def}
Let \( K \subset \Sigma \) compact. We call a sequence of relatively open sets \( \{U_k\}_{k=1}^{\infty} \) in \( \Sigma \) an in-approximation of \( K \) if
    \begin{itemize}
    \item \( U_k \subset U_{k+1}^{rc} \) for all \( i \);
    \item \( \sup_{X \in U_k} \text{dist}(X, K) \rightarrow 0 \text{ as } k \rightarrow \infty \).
\end{itemize}
\end{definition}

\begin{theorem}[\protect{\cite{MR1728376} and \cite[Theorem 1.4]{MR3740399} and \cite[Proposition 3.4 and Theorem 3.5]{hab_thesis}}]\label{in_approx_thm}
Assume \( \Omega \subset \mathbb{R}^2 \) is a bounded domain. Let \( K \subset \Sigma \) be compact and assume that \( \{U_k\}_{k=1}^{\infty} \) is an in-approximation of \( K \) (see \Cref{in_approx_def}). Then for each piecewise affine Lipschitz map \( v\colon \Omega \rightarrow \mathbb{R}^2 \) with \( Dv(x) \in U_1 \) in \( \Omega \) there exists a Lipschitz map \( u\colon \Omega \rightarrow \mathbb{R}^2 \) verifying
\[
    Du(x) \in K \quad \text{a.e. in } \Omega, \quad u(x) = v(x) \text{ on } \partial\Omega.
\]  
\end{theorem}

\subsection{A computational approach}
\hspace{.3in}

We also utilize an additional, earlier, theory of 
Sz\'{e}kelyhidi \cite{MR2048569} where he develops the following idea: one method of constructing convex integration solutions to a particular class of systems would be to fix the system under consideration (such as the $p$-system \eqref{psystem} with a \emph{fixed} function $p$), and then look for $N$ points in the corresponding constitutive set which are in $T_N$ configuration. Another, and more practical approach would be to computationally search for a $T_N$ configuration which lives in the constitutive set for some system (in the appropriate class) which is constructed \emph{ad-hoc}. In our setting (i.e. the system \eqref{psystem}), this would correspond to finding a large $T_5$, and then only afterwards determining the function $p\colon \mathbb{R}\to\mathbb{R}$. 

\subsection{Plan for the paper}

The sketch of our proof then is the following. First, we computationally find a large $T_5$ which has algebraic constraints which allow for it to be contained in the constitutive set $\mathcal{K}$ for at least one genuinely nonlinear, strictly hyperbolic $p$-system \eqref{psystem}. We then construct an appropriate ``subsolution'' which, via the convex integration and in-approximation theory, will give rise to the ``full'' solutions with Riemann initial data and verifying the conclusions of the Main Theorem (\Cref{main_theorem}).

In \Cref{sec:sol_by_convex_int} we introduce rank-one convexity and $T_N$ configurations. In \Cref{sec:rank_one_shocks} we introduce shocks and admissibility conditions on shocks. In \Cref{sec:proof_main_thm} we give the proof of the Main Theorem (\Cref{main_theorem}), which will follow from the Main ``MATLAB'' Proposition (\Cref{matlab_lemma}), proven in \Cref{sec:T_N_existence_section}.

\section{Solutions via convex integration}\label{sec:sol_by_convex_int}

We follow the framework from \cite{MR2008346}: For a given function $p\colon\mathbb{R}\to\mathbb{R}$, we consider stream functions \( \psi(x, t) \colon \mathbb{R}^2 \rightarrow \mathbb{R}^2 \) such that
\begin{equation}
\begin{aligned}\label{correspondance_psi}
(v, -u ) &= ((\psi_1)_x, (\psi_1)_t) \\
(u, -p(v)) &= ((\psi_2)_x, (\psi_2)_t),
\end{aligned}
\end{equation}
for functions $v,u\colon\mathbb{R}\to\mathbb{R}$ and where we write \( \psi = (\psi_1, \psi_2) \) for the components of \( \psi \). In terms of \( \psi \), the  system \eqref{psystem} is equivalent to the first order differential inclusion 
\begin{equation}\label{diff_inclusion}
D\psi \in \mathcal{K},
\end{equation}
where the constitutive set $\mathcal{K}\subset\mathbb{R}^{2\times2}$ is given by
\begin{equation}\label{K_def}
\mathcal{K} := \left\{ G(v,u)
: v,u \in \mathbb{R} \right\},
\end{equation}
where 
\begin{align}\label{G_def}
    G(v,u)\coloneqq \begin{pmatrix}
v & -u \\
u & -p(v)
\end{pmatrix}.
\end{align}

In this paper, we will construct such stream functions \( \psi \) (via \Cref{in_approx_thm}). They will be only Lipschitz, so their derivatives, which are solutions to \eqref{psystem} via the correspondence \eqref{correspondance_psi}, will be only $L^\infty$. 

\subsection{$T_N$ configurations and rank-one convexity}\label{sec:TN_rank_one}

 In this section we revisit the important definitions and results on \emph{rank-one convexity}. A function $f\colon\mathbb{R}^{m\times n}\to\mathbb{R}$ is rank-one convex if $f$ is convex along each rank-one line. The \emph{rank-one convex hull} of a set of matrices is defined by separation with rank-one convex functions: For a compact set $K\subset \mathbb{R}^{m \times n}$, we define the rank-one convex hull
\begin{align}
K^{\text{rc}}\coloneqq \big\{X\in\mathbb{R}^{m\times n} : f(X)\leq \sup_{K} f \mbox{ for all } f\colon \mathbb{R}^{m\times n}\to\mathbb{R} \mbox{ rank-one convex}\big\}.
\end{align}

Let us denote by \( \{X_1, \ldots, X_N\} \) the unordered set of matrices \( X_i, i = 1, \ldots, N \) and by \( (X_1, \ldots, X_N) \) the ordered \( N \)-tuple.

\begin{definition}[\( T_N \)-configuration)]\label{TN_Def}
    Let \( X_1, \ldots, X_N \in \mathbb{R}^{m \times n} \) be \( N \) matrices such that $\mathrm{rank}(X_i - X_j) > 1$ for all \( i \neq j \). The ordered set \( (X_1, \ldots, X_N) \) is said to be in \( T_N \) configuration if there exist \( P, C_i \in \mathbb{R}^{m \times n} \) and \( \kappa_i > 1 \) such that

\begin{equation}
    \begin{aligned}\label{TN_param}
X_1 &= P + \kappa_1 C_1 \\
X_2 &= P + C_1 + \kappa_2 C_2 \\
& \vdots \\
X_N &= P + C_1 + \ldots + C_{N-1} + \kappa_N C_N,
\end{aligned}
\end{equation}
and furthermore \( \mathrm{rank} \, C_i \) = 1 and \( \sum_{i=1}^{N} C_i = 0 \).
\end{definition}  

The reason $T_N$ configurations are significant is because their cyclic structure increases the size of the rank-one convex hull, as made precise by the following lemma which is well known in the literature (see e.g. \cite{MR1983780,MR1320538}).

\begin{figure}[tb]
      \includegraphics[width=.7\textwidth]{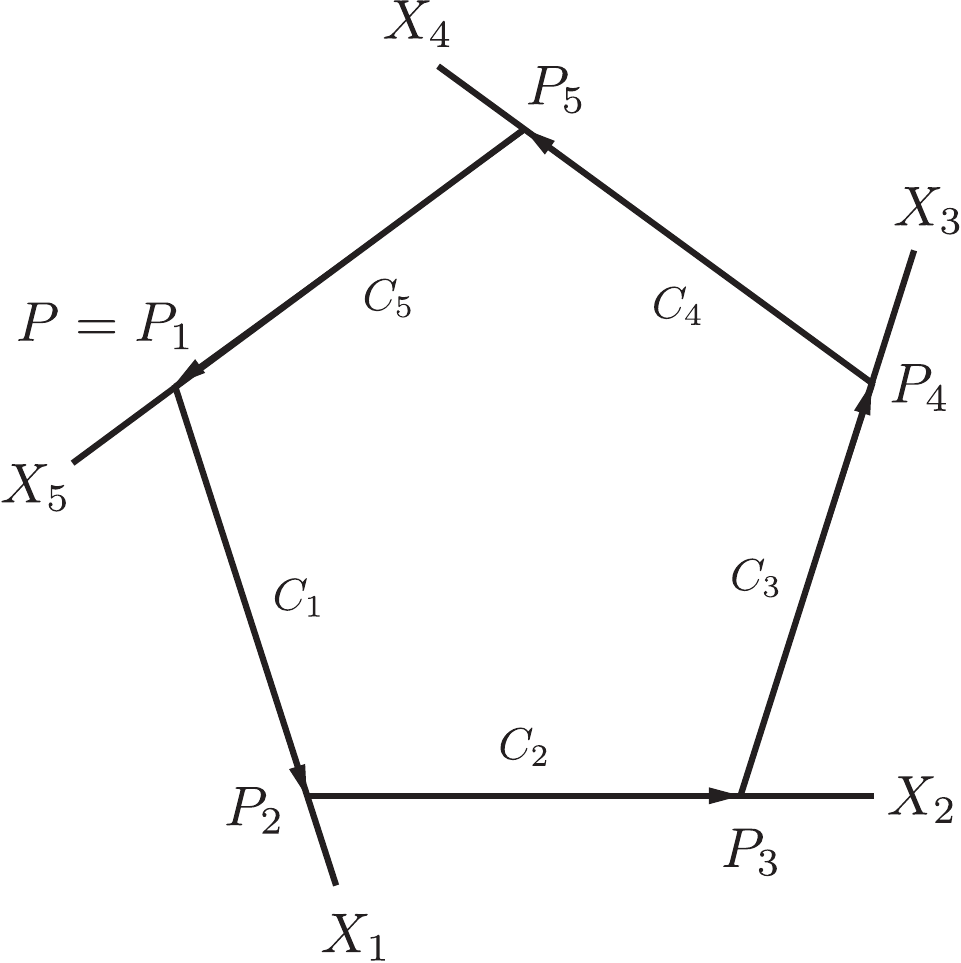}\hspace{.8in}
  \caption{A schematic of a $T_5$ configuration.}\label{fig:example_T5}
\end{figure}

\begin{lemma}\label{convex_hull_TN}
Assume \( (X_i)_{i=1}^N \) is a \( T_N \)-configuration. Then
\[
\{P_1, \ldots, P_N\} \subset \{X_1, \ldots, X_N\}^{rc},
\]
where \( P_1 \coloneqq P \) and \( P_i \coloneqq P + \sum_{j=1}^{i-1} \lambda_j C_j \text{ for } i = 2, \ldots, N \).
\end{lemma} 

See \Cref{fig:example_T5} for a diagrammatic view of a $T_N$, for $N=5$.

\Cref{TN_Def} does not give any insight into whether or not a particular (ordered) set of $N$ matrices are in $T_N$ configuration or not. However,  we have the following powerful characterization from \cite{MR2118899}:

\begin{theorem}[Algebraic criterion \cite{MR2118899}]\label{alg_crit_theorem}
For any \( \mu \in \mathbb{R} \) and any \( A \in \mathbb{R}^{N\times N}_{\text{sym}} \) with \( A_{ii} = 0 \) for \( i = 1, \ldots, N \), we first define
\[
A^{\mu} \coloneqq \begin{pmatrix}
0 & A_{12} & A_{13} & \cdots & A_{1N} \\
\mu A_{12} & 0 & A_{23} & \cdots & A_{2N} \\
\vdots & \vdots & \vdots & \ddots & \vdots \\
\mu A_{1N} & \mu A_{2N} & \mu A_{3N} & \cdots & 0
\end{pmatrix}.
\]

Then, let \( (X_1, \ldots, X_N) \in (\mathbb{R}^{2\times2})^N \) and let \( A \in \mathbb{R}^{N\times N} \) with \( A_{ij} = \det(X_i - X_j) \). Then \( (X_1, \ldots, X_N) \) are in \( T_N \)-configuration if and only if there exist \( \lambda_1, \ldots, \lambda_N > 0 \) and \( \mu > 1 \) such that \( A^{\mu}\lambda = 0 \), where $\lambda\in\mathbb{R}^N$ is defined as $\lambda\coloneqq(\lambda_1, \ldots, \lambda_N)$.
\end{theorem}

In fact, from \( \mu \) and \( \lambda = (\lambda_1, \ldots, \lambda_N) \) we can easily compute the parametrization
\( (P, C_i, \kappa_i) \) (see \eqref{TN_param}) of the \( T_N \)-configuration \( (X_1, \ldots, X_N) \). In detail, using the definition of the \( P_i \) from \Cref{convex_hull_TN}, we have (from \cite{MR2118899}):

\begin{equation}
\begin{aligned}\label{compute_param}
P_1 &= \frac{1}{\lambda_1 + \cdots + \lambda_N} (\lambda_1 X_1 + \ldots + \lambda_N X_N) \\
P_2 &= \frac{1}{\mu\lambda_1 + \lambda_2 + \cdots + \lambda_N} (\mu\lambda_1 X_1 + \lambda_2 X_2 + \cdots + \lambda_N X_N) \\
& \vdots \\
P_N &= \frac{1}{\mu\lambda_1 + \cdots + \mu\lambda_{N-1} + \lambda_N} (\mu\lambda_1 X_1 + \cdots + \mu\lambda_{N-1} X_{N-1} + \lambda_N X_N)
\end{aligned} 
\end{equation}

\subsection{Stability of $T_5$ configurations in \( \mathbb{R}^{2\times2} \) }

We have the following stability result for $T_5$ configurations in \( \mathbb{R}^{2\times2} \)

\begin{lemma}[\protect{Stability of $T_5$ configurations in \( \mathbb{R}^{2\times2} \) \cite[Lemma 2.4]{MR3740399}}]\label{T5_stab} Let \( (X_1, \ldots, X_5) \) be a \( T_5 \)-configuration in \( \mathbb{R}^{2\times2} \) with \( \det(X_i - X_j) \neq 0 \) for all \( i \neq j \). Then there exists \( \epsilon > 0 \) so that any \( (\widetilde{X}_1, \ldots, \widetilde{X}_5) \) with \( |\widetilde{X}_i - X_i| < \epsilon, i = 1 \ldots 5 \), is also in \( T_5 \)-configuration.
\end{lemma}

\subsection{Large $T_5$ configurations}
We can now introduce the main definition from \cite{MR3740399}, which we will utilize in our paper.

It should be noted that a set consisting of five matrices has the potential to generate multiple $T_5$-configurations, with each configuration arising from a distinct permutation of the set's elements. 

In order to study these scenarios, let \( \{X_1^0, \ldots, X_5^0\} \) be a 5-element set and let \( S_5 \) be the permutation group of 5 elements. To any \( \sigma \in S_5 \) is associated a 5-tuple \( (X_{\sigma(1)}^0, \ldots, X_{\sigma(5)}^0) \). If this 5-tuple is a \( T_5 \)-configuration, then there exists a map
\[
  (X_{\sigma(1)}^0, \ldots, X_{\sigma(5)}^0) \mapsto (P^{\sigma}_{\sigma(1)}, \ldots, P^{\sigma}_{\sigma(5)})
\]
where \( P^{\sigma}_{\sigma(i)} \) are the corresponding matrices from \Cref{convex_hull_TN}, so that in particular
\[
  \text{rank}(P^{\sigma}_{\sigma(i)} - X_{\sigma(i)}^0) = 1 \quad \text{and} \quad P^{\sigma}_{\sigma(i)} \in \{X_1^0, \ldots, X_5^0\}^{rc}.
\]
Let
\[
  C_i^{\sigma} := P^{\sigma}_{\sigma(i)} - X_i^0.
\]
See \cite[p.~19]{MR3740399} for more details.

\begin{definition}[Large $T_5$]\label{def:largeT5}
    We call a 5-element set \( \{X_1^0, \ldots, X_5^0\} \subset \mathbb{R}^{2 \times 2} \) a \emph{large \( T_5 \)-set} if there exist at least three permutations \( \sigma_1, \sigma_2, \sigma_3 \) such that \( (X_{\sigma_j(1)}^0, \ldots, X_{\sigma_j(5)}^0) \) is a \( T_5 \)-configuration for each \( j = 1, 2, 3 \), and furthermore the corresponding rank-one matrices \( C_i^{\sigma_1}, C_i^{\sigma_2}, C_i^{\sigma_3} \) (see \eqref{TN_param}) are linearly independent for all \( i = 1, \ldots, 5 \).
\end{definition} 

Considering the stability presented in \Cref{T5_stab}, it is clear that large $T_5$  sets maintain stability when subject to minor perturbations.

We can now present the major result from \cite{MR3740399}.

\begin{theorem}[\protect{Existence of in-approximations for large $T_5$s configurations \cite[Theorem 2.8]{MR3740399}}]\label{thm:T5_large_set_in_approx}
Let \( K = \{X_1^0, \ldots, X_5^0\} \) be a large \( T_5 \) set. Then there exists an in-approximation \( \{U_k\}_{k=1}^{\infty} \) of \( K \).
\end{theorem}

\section{Rank-one connections and shocks}\label{sec:rank_one_shocks}

For  $(v_L,u_L)$ and $(v_R,u_R) \in \mathbb{R}^{2}$ and \( \sigma \in \mathbb{R} \), the function
\begin{equation}\label{shock_sol}
(v(x,t),u(x, t)) \coloneqq \begin{cases}
(v_L,u_L) & \text{if } x < \sigma t, \\
(v_R,u_R) & \text{if } x > \sigma t,
\end{cases}
\end{equation}
is a weak solution to \eqref{psystem} if and only if $(v_L,u_L)$, $(v_R,u_R)$ and \( \sigma \) verify the Rankine-Hugoniot jump condition,
\begin{equation}
\begin{aligned}\label{RH_cond}
    u_R-u_L&=\sigma(v_R-v_L)\\
    p(v_R) - p(v_L) &= \sigma(u_R - u_L).
\end{aligned}
\end{equation}
In this case, \eqref{shock_sol} is called a \emph{shock solution} to \eqref{psystem}. The scalar $\sigma$ is the \emph{shock speed}. In terms of the set $\mathcal{K}$ (see \eqref{K_def} and \eqref{G_def}), $(v_L,u_L)$, $(v_R,u_R)$ are a shock if and only if
\begin{align}\label{rank_one_shock}
    \det(G(v_R,u_R)-G(v_L,u_L))=0,
\end{align}
i.e. $G(v_R,u_R)$ and $G(v_L,u_L)$ are rank-one connected.

\subsection{Hugoniot locus}
From the Rankine-Hugoniot condition, we can define the \emph{Hugoniot locus} at a point $(v_L,u_L)$ as follows,
\begin{align}\label{locus}
    H(v_L,u_L)\coloneqq \big\{(v_R,u_R)\in\mathbb{R}^2 \, | \, f(v_R,u_R)-f(v_L,u_L)=\sigma( (v_L,u_L)-(v_R,u_R)) \\ 
    \hspace{2in}\mbox{ for some }\sigma\in\mathbb{R}\big\},
\end{align}
where the flux $f$ is defined as $f(v,u)=(u,p(v))$ in our case for the $p$-system \eqref{psystem}.

Locally around $(v_L,u_L)$ the Hugoniot locus $H(v_L,u_L)$ is the union of two smooth curves (the same holds for general strictly hyperbolic $n\times n$ systems, where the Hugoniot locus is the union of $n$ smooth curves -- see e.g. \cite[Section 8.2]{dafermos_big_book}). Let us label the two curves \( S^1_{(v_L,u_L)} \) and \( S^2_{(v_L,u_L)} \).  Each curve passes through the point \( (v_L,u_L) \). These are defined as the \emph{shock curves} of the Hugoniot locus. Let us (smoothly) parameterize \( S^1_{(v_L,u_L)} \) and \( S^2_{(v_L,u_L)} \) as follows: \( S^1_{(v_L,u_L)} = S^1_{(v_L,u_L)}(s) \) and \( S^2_{(v_L,u_L)} = S^2_{(v_L,u_L)}(s) \) with \( S^1_{(v_L,u_L)}(0) = S^2_{(v_L,u_L)}(0) = (v_L,u_L) \). We also choose a (smooth) parameterization for the speed \( \sigma \). To be precise, we choose \( \sigma^1_{(v_L,u_L)} \) and \( \sigma^2_{(v_L,u_L)} \) such that
\begin{equation}
\sigma^k_{(v_L,u_L)}(s)(S^k_{(v_L,u_L)}(s) - (v_L,u_L)) = f(S^k_{(v_L,u_L)}(s)) - f(v_L,u_L),
\end{equation}
for \( k = 1,2 \).

\subsection{Lax \emph{E}-condition and Liu \emph{E}-condition}\label{sec:shock_conds}

For a general system, a shock
\begin{align}
    (U_L,S^k_{U_L}(s_R),\sigma^k_{U_L}(s_R))
\end{align}
(from characteristic field $k$) is said to verify the \emph{Lax \emph{E}-condition} if
\begin{align}\label{lax_cond}
    \lambda_k(S^k_{U_L}(s_R))\leq \sigma^k_{U_L}(s_R) \leq \lambda_k(U_L).
\end{align}
See \cite[p.~274]{dafermos_big_book} for more on the Lax \emph{E}-condition.

The shock $(U_L,S^k_{U_L}(s_R),\sigma^k_{U_L}(s_R))$ is said to verify the \emph{Liu \emph{E}-condition} if
\begin{align}\label{liu_cond}
    \sigma^k_{U_L}(s_R) \leq  \sigma^k_{U_L}(s), \text{ for all $s$ between $0$ and $s_R$}.
\end{align}

The Liu \emph{E}-condition is in some sense stricter and more discriminating than the Lax \emph{E}-condition, in particular for systems which are not genuinely nonlinear. See \cite[Section 8.4]{dafermos_big_book} for more on the Liu \emph{E}-condition.

The $p$-system \eqref{psystem} admits two families of shocks which both verify the Lax \emph{E}-condition and the Liu \emph{E}-condition, in the case when the $p$-system is genuinely nonlinear, i.e. $p''\neq0$. For a reference, see \cite[p.~275 and p.~280]{dafermos_big_book}.

\section{Proof of Main Theorem (\Cref{main_theorem})}\label{sec:proof_main_thm}

The proof of the Main Theorem (\Cref{main_theorem}) will follow from the following Main Proposition and an ancillary Lemma.

\begin{proposition}[Main Proposition -- The MATLAB Proposition]\label{matlab_lemma}
There exists a large $T_5$ configuration (see \Cref{def:largeT5}) $\{X_1^0,\ldots,X_5^0\}\subset \mathbb{R}^{2\times2}$ and scalars $D_1,\ldots,D_5<0$  verifying the inequalities 

\begin{align}\label{inequalities_for_convexity_extension_lemma}
(X^0_j)_{2,2}-(X^0_i)_{2,2} > D_i\cdot \Big((X^0_j)_{1,1}-(X^0_i)_{1,1}\Big) \mbox{ for all } i\neq j,
\end{align}

\begin{align}\label{symmetry_T5}
(X^0_i)_{2,1}+(X^0_i)_{1,2}=0  \mbox{ for all } i=1,\ldots,5.
\end{align}

Furthermore, the large $T_5$ configuration contains no rank-one connections, i.e. $\rm{rank}(X_i^0-X_j^0)=2$ for all $i\neq j$.

A concrete example of a large $T_5$ is in the appendix (\Cref{sec:T5_example}).
\end{proposition}
\begin{remark}
In the work \cite{MR2048569}, Sz\'{e}kelyhidi introduces a computer-assisted proof technique for finding systems which admit $T_N$ configurations.
    
    In \cite{MR2048569}, the existence of the required $T_5$ configuration was reduced to solving a linear system of inequalities -- using the \emph{simplex algorithm} in Maple V. To get a \emph{linear} system of inequalities, Sz\'{e}kelyhidi writes the $C_i$ (as in the parameterization of a $T_N$ configuration, see \eqref{TN_param}), as the tensor product of two vectors $C_i=\alpha_i\bigotimes\beta_i$. However, in order for the constraint $\sum_i C_i=0$ to be linear in the entries of the vectors $\alpha_i$ and $\beta_i$, one of the vectors must be kept fixed. This reduces the flexibility on the $T_N$ configuration and makes it more difficult to impose additional algebraic constraints on the $T_N$ configuration (such as  \eqref{inequalities_for_convexity_extension_lemma} and \eqref{symmetry_T5}). See \cite[p.~145]{MR2048569} for details. 

     In this paper, we expand on the technique of \cite{MR2048569} and introduce a novel method of finding $T_N$ configurations with \emph{nonlinear algebraic constraints}. 
\end{remark}

The proof and discussion of the MATLAB code is in \Cref{sec:T_N_existence_section}, below.

We will also make use of the following fact about convex functions (cf. \eqref{inequalities_for_convexity_extension_lemma}, above):

\begin{lemma}[Algebraic inequalities which yield a strictly convex function]\label{eta_extension_prop}
Fix $n,N\in\mathbb{N}$. Assume there exists $(x_i,h_i,D_i)\in\mathbb{R}^n\times\mathbb{R}\times\mathbb{R}^{n}$, for $i=1,\ldots,N$, such that the strict inequalities
\begin{align}\label{inequalities_for_convexity}
h_j > h_i + D_i\cdot (x_j-x_i) \mbox{ for all } i\neq j
\end{align}
are verified. Then, there exists a smooth and strictly convex function $\eta\colon \mathbb{R}^n\to\mathbb{R}$ such that $\eta(x_i)=h_i$ and $D\eta(x_i)=D_i$, for all $i$.
\end{lemma}
\begin{remark}
It is well known that such an $\eta$ exists if we only require it be convex. For instance, see \cite[p.~143]{MR2048569}. However, in this \Cref{eta_extension_prop} $\eta$ is \emph{strictly} convex. 
\end{remark}

The proof of \Cref{eta_extension_prop} is in \Cref{sec:eta_extension_prop}, below.

We can now begin the proof of \Cref{main_theorem}.

\subsection{Proof of \Cref{main_theorem}}

Consider the large $T_5$ $\{X_1^0,\ldots,X_5^0\}$ from \Cref{matlab_lemma}. Consider the particular ordering $(X_1^0,\ldots,X_5^0)$, as well as the parameterization $(P,C_i,\kappa_i)$ associated to this ordering (see \eqref{TN_param}).

\uline{Step 1}  From \eqref{inequalities_for_convexity_extension_lemma} and \eqref{symmetry_T5}, as well as \Cref{eta_extension_prop}, we have a smooth function $p\colon \mathbb{R}\to\mathbb{R}$ verifying $p'>0$ and $p''<0$ and such that under the identification \eqref{correspondance_psi}, the $T_5$ given by $\{X_1^0,\ldots,X_5^0\}$ lives in the the constitutive set $\mathcal{K}$ (see \eqref{K_def}).

\uline{Step 2}

Remark that $P-X_1^0$ and $P-X_5^0$ are both rank-one matrices. Thus it is possible for a Lipschitz function to have a derivative taking only the values $P$, $X_1^0$, and $X_5^0$.

\begin{figure}[tb]
      \includegraphics[width=.6\textwidth]{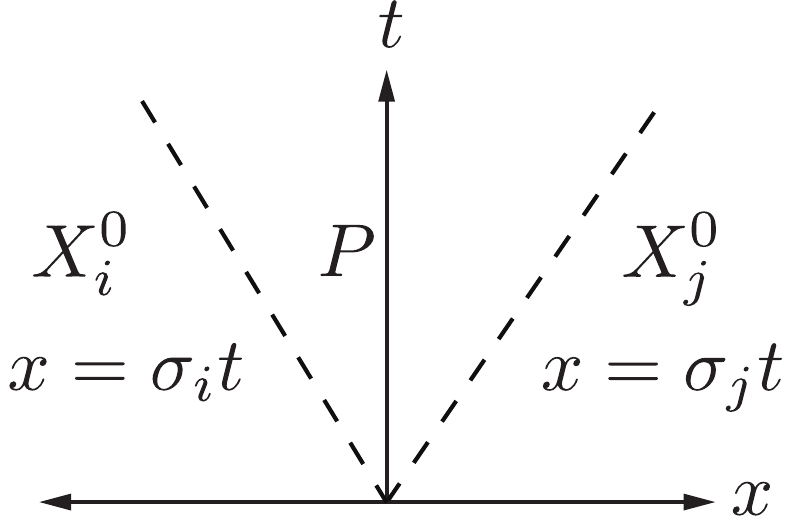}
  \caption{A diagram of values of $DS$, the Jacobian of the subsolution $S$.}\label{fig:wedge}
\end{figure}

In particular, define a ``wedge'' Lipschitz function $S(x,t)\colon \mathbb{R}\times[0,\infty)\to\mathbb{R}^{2}$, such that the derivative of $S$ verifies
\begin{align}\label{def_S}
   DS(x,t) = 
  \begin{cases} 
   X_i^0 & \text{if } x < \sigma_i t, \\
   P       & \text{if } \sigma_i t < x < \sigma_j t,\\
   X_j^0 & \text{if } \sigma_j t < x,
  \end{cases} 
\end{align}
where $DS$ is the Jacobian of $S$ in the $x$ and $t$ variables, $\sigma_i$ is the value of the slope needed in the $(x,t)$ plane to allow for a Lipschitz function taking adjacent derivative values $P$, $X_i^0$ (as in \eqref{def_S}) and likewise $\sigma_j$ is the value of the slope needed to allow for a Lipschitz function taking adjacent derivative values $P$, $X_j^0$. See \Cref{fig:wedge}. Our use of the ``wedge'' \eqref{def_S} is closely linked with the \emph{fan subsolution} used in \cite[Definition 3.4]{multi_d_illposed}.

We have two choices for the $i$ and $j$. Either $i=1,j=5$ or $i=5,j=1$. We must choose such that $\sigma_i<\sigma_j$. It is possible to check directly the large $T_5$ we found computationally (see \Cref{sec:T5_example}), and indeed for this particular large $T_5$ we find $\sigma_1<\sigma_5$.

\uline{Step 3}

From \Cref{thm:T5_large_set_in_approx}, we get an in-approximation $\{U_k\}_{k=1}^{\infty}$ for the large $T_5$ $\{X_1^0,\ldots,X_5^0\}$. Remark from the proof of \Cref{thm:T5_large_set_in_approx} in  \cite{MR3740399} that it is possible to choose the 
in-approximation $\{U_k\}_{k=1}^{\infty}$ such that $P\in U_1$ (see \cite[p.~19]{MR3740399} and also \Cref{convex_hull_TN}).

We then apply \Cref{in_approx_thm}. Remark that we can apply \Cref{in_approx_thm} because after multiplying the second column of each matrix in our large $T_5$ by $(-1)$, the large $T_5$ will live in the set of symmetric $2\times 2$ matrices. Remark also that after  multiplying the second column of each matrix in our large $T_5$ by $(-1)$, we still have a $T_5$ structure because if $C_i$ is a rank-one matrix from within the context of \Cref{TN_Def}, then after multiplying one of its columns by $(-1)$, it will still be rank-one.  This gives us a Lipschitz map $\mathbb{U}(x,t)\colon \mathbb{R}\times[0,\infty)\to\mathbb{R}^{2\times 2}$.
Notice that, in the language of convex integration, $S$ is the ``subsolution'' and $\mathbb{U}$ is the full solution to the differential inclusion \eqref{diff_inclusion}.
Remark that under the identification \eqref{correspondance_psi}, with $v\coloneqq \mathbb{U}_{1,1}$ and $u\coloneqq \mathbb{U}_{2,1}$, the 2-tuple $(v,u)$ is a solution to \eqref{psystem}. Note the solution $(v,u)$ is uniformly bounded. This follows because $\mathbb{U}$ is Lipschitz.

Due to the wedge shape \eqref{def_S}, it follows that $(v(\cdot,t),u(\cdot,t))$ converges in $L^1_{\text{loc}}$ to the Riemann initial data
\begin{align}
(v^0,u^0)=
\begin{cases}
       ((X_i^0)_{1,1},(X_i^0)_{2,1}) & \text{if } x < 0, \\
   ((X_j^0)_{1,1},(X_j^0)_{2,1}) & \text{if }  x > 0,
   \end{cases}
\end{align}
as $t\to0^{+}$. Thus, $(v,u)$ is a weak solution to \eqref{psystem} with Riemann initial data.

\uline{Step 4} Recall from \Cref{matlab_lemma} that the $T_5$ configuration $\{X_1^0,\ldots,X_5^0\}$ contains no rank-one connections, i.e. $\rm{rank}(X_i^0-X_j^0)=2$ for all $i\neq j$. Thus, recalling that a rank-one connection corresponds to a shock solution (see \Cref{sec:rank_one_shocks}), we conclude that the solution $(v,u)$ contains no shocks, and thus vacuously meets the criteria for the Lax \emph{E}-condition (see \eqref{lax_cond}) and Liu \emph{E}-condition (see \eqref{liu_cond}).

\uline{Step 5} Remark that due to \Cref{T5_stab}, after we have fixed the function $p$, we can perturb the large $T_5$ $\{X_1^0,\ldots,X_5^0\}$ to find different large $T_5$ configurations which also live in the constitutive set $\mathcal{K}$. In doing so, we can perturb the initial data \eqref{riemann_initial} within an open set. Moreover, for each fixed Riemann initial data, we can perturb the \emph{other three} values of our solution (not used in the initial data), thus yielding infinitely many solutions for each fixed Riemann problem.

Finally, from the proof of \Cref{in_approx_thm}, it is clear that for each fixed time $t>0$, within the wedge $\sigma_i t < x < \sigma_j t$, the solution $(v(\cdot,t),u(\cdot,t))$ will have infinite total variation.

This concludes the proof of the Main Theorem.

\subsection{Proof of Main Proposition -- \Cref{matlab_lemma} (existence of large $T_5$)}\label{sec:T_N_existence_section}

In this section, we show the existence of a suitable large $T_5$ configuration. In particular, we prove \Cref{matlab_lemma}.


Our proof is computer-assisted; the MATLAB code is on the GitHub\footnote{See \url{https://github.com/sammykrupa/BV-blowup-for-p-system}}.

\uline{Step 1} We first enter into MATLAB the parameterization of a $T_N$ configuration 
\begin{align}
    \{X_1^0,\ldots,X_5^0\},
\end{align} 
in the ordering $(X_1^0,\ldots,X_5^0)$ (as in \eqref{TN_param}). In particular, the rank-one matrices $C_i$ (from \eqref{TN_param}) are written as the tensor products $C_i=a_i\bigotimes n_i$, for $a_i\in\mathbb{R}^{2}$, $n_i\in\mathbb{R}^{2}$ ($i=1,\ldots,5$). From the definition of a $T_N$ configuration, the constraint $\sum C_i =0$ is also input into MATLAB. Further, we want the $X_i^0$ (playing the role of the $X_i$ in the definition of $T_N$ configuration) to satisfy  \eqref{inequalities_for_convexity_extension_lemma} and \eqref{symmetry_T5}.  

For simplicity, we choose fixed values for the $\kappa_i$ and we also assume that $P=0$ (where $P$ is from \eqref{TN_param}) and thus the solver will not have to determine these values. 

Then, the MATLAB R2023b solver \emph{fmincon} and the \emph{interior-point} algorithm (see \cite{MATLAB_solve,fmincon}) return \emph{numeric} values for the $a_i$, $n_i$, and $D_i$ in double precision (see e.g. \cite{MATLAB_double}). These values of $a_i$, $n_i$, and $D_i$, along with the fixed choices for the $\kappa_i$, are then converted to exact \emph{symbolic values} within MATLAB (see \cite{MATLAB_symbolic}). Label these symbolic values $\hat a_i$, $\hat n_i$, $\hat \kappa_i$, and $\hat D_i$.

Remark that solving a nonlinear system of inequalities using MATLAB with fmincon and the interior-point algorithm,  gives different solutions depending on the initial point at which the solver starts. Our MATLAB code provides a fixed initial point which returns the large $T_5$ example which we use in this paper (see the appendix, \Cref{sec:T5_example}). The code can also choose a random initial point for the solver each time the code is run. Not all initial points lead to feasible solutions.

From this point on, we perform all further computations symbolically within MATLAB (see \cite{MATLAB_symbolic}). This allows for rigorous mathematical statements.  

\uline{Step 2} We then calculate the points $X_1^0,\ldots, X_5^0$ in a $T_5$ configuration by using the $\hat a_i$, $\hat n_i$, $\hat \kappa_i$ in the formulas for the parameterization of a $T_5$ (see \eqref{TN_param}). As above, continue to assume that $P=0$.  We then check symbolically that the constraints \eqref{inequalities_for_convexity_extension_lemma} and \eqref{symmetry_T5} hold, as well as $\rm{rank}(X_i^0-X_j^0)=2$ for all $i\neq j$.

\uline{Step 3} Next, it is necessary to determine which orderings of the matrices $\{X_1^0,\ldots,X_5^0\}$ give a $T_5$ configuration. To do so, the MATLAB code uses \Cref{alg_crit_theorem}\footnote{To determine which orderings we should test symbolically with \Cref{alg_crit_theorem}, we first used trial and error and checked many orderings at random.}. This involves symbolically solving for the roots of the polynomial (in $\mu$)
\begin{align}
    \det(A^\mu)=0,
\end{align}
where $A^\mu$ is from the context of \Cref{alg_crit_theorem}.

Remark that due the approximate nature of the MATLAB solver fmincon, MATLAB does \emph{not} return values for $\hat a_i$ and $\hat n_i$ which make a rank-one polygon. In other words, $\sum_i \hat a_i \bigotimes \hat n_i \neq 0$ -- there is indeed a very small numerical error. Thus even for the ordering of the $\{X_1^0,\ldots,X_5^0\}$ which the numerical solver found (i.e. $(X_1^0,\ldots,X_5^0)$), symbolically applying \Cref{alg_crit_theorem} will give a slightly different parameterization (including possibly with $P\neq0$).

\uline{Step 4} From applying \Cref{alg_crit_theorem} to each ordering, one can compute the parameterization $(P,C_i,\kappa_i)$ of each ordering (see \eqref{compute_param}). From this, it is necessary to check the linear independence of the rank-one directions $C_i$ needed for a large $T_5$ (see \Cref{def:largeT5}). This completes the proof.

Steps 1-4 of the proof are implemented in the symbolic MATLAB code available on the GitHub.

\subsection{Proof of \Cref{eta_extension_prop} (existence of a strictly convex function)}\label{sec:eta_extension_prop}

We use ideas from the proof for the convex case given in \cite[p.~143]{MR2048569} (but the function in \cite[p.~143]{MR2048569} is not necessarily \emph{strictly} convex).

We define $\eta_0\colon\mathbb{R}^n\to\mathbb{R}$, as follows
\begin{align}\label{max_for_convexity}
\eta_0(x)\coloneqq \max_{i}\{h_i+D_i\cdot(x-x_i)+\epsilon_0\abs{x-x_i}^2\},
\end{align}
for a small $\epsilon_0>0$. Due to the strictness in the finitely many inequalities \eqref{inequalities_for_convexity}, we can choose an $\epsilon_0$ sufficiently small such that $\eta_0(x_i)=h_i$ for all $i$.

Consider $m$, a smooth mollifier on $\mathbb{R}^n$, positive, supported on a small ball around the origin, verifying $\int m(x)\,dx =1$ and  verifying $\int x m(x)\,dx =0$. 

Then, due again to the strictness in the finitely many inequalities \eqref{inequalities_for_convexity}, we can choose the support of $m$ sufficiently small, and possibly reduce $\epsilon_0$ further, such that locally around $x_i$, we have 
\begin{align}
m\ast \eta_0(x)&=\int \Big(h_i+D_i\cdot[(x-y)-x_i]+\epsilon_0\abs{(x-y)-x_i}^2\Big)m(y)\,dy\\
&=h_i+D_i\cdot(x-x_i) + \int \epsilon_0\abs{(x-y)-x_i}^2 m(y)\,dy
\end{align}
In particular,
\begin{align}
m\ast \eta_0(x_i)=h_i + \int \epsilon_0\abs{y}^2 m(y)\,dy.
\end{align}

We want $m\ast \eta_0(x_i)=h_i$. The term $ \int \epsilon_0\abs{y}^2 m(y)\,dy$ is error. However, once again due to the strictness in the inequalities \eqref{inequalities_for_convexity}, and the fact that there are only finitely many inequalities, we can modify the $h_i$, replacing them with slightly smaller values. In particular, for  $\tilde{h}_i$ sufficiently close to $h_i$ and verifying $\tilde{h}_i<h_i$, the inequalities \eqref{inequalities_for_convexity} are satisfied with the $h_i$ replaced by $\tilde{h}_i$. This would then give,
\begin{align}
m\ast \eta_0(x_i)=\tilde{h}_i + \int \epsilon_0\abs{y}^2 m(y)\,dy.
\end{align}

The quantity $\int \epsilon_0\abs{y}^2 m(y)\,dy$ goes to zero as the support of the mollifier becomes smaller. Thus, we can choose $\tilde{h}_i$ sufficiently close to $h_i$ and a support of the mollifier $m$ sufficiently small, such that $m\ast \eta_0(x_i)=h_i$. Note also that due to the map $\mathbb{R}^n\ni x\mapsto \norm{x}^2$ having zero first derivative at the origin, $D (m\ast \eta_0)(x_i)=D_i$.

We then take as our definition, 
\begin{align}
\eta\coloneqq m\ast \eta_0.
\end{align}

To conclude the proof of the lemma we show that $\eta$ is strictly convex. To do this, we show that for each $x_0\in\mathbb{R}^n$ and $v\in\mathbb{M}^{n}$, the map $\mathbb{R}\ni t\mapsto \eta(x_0+tv)$ has positive second derivative at $t=0$. 

Let us comment first on why this will show that $\eta$ is strictly convex.

We calculate,
\begin{align}
\frac{d^2}{dt^2} \eta(x_0+tv) =v^T D^2 \eta (x_0+tv) v.
\end{align}

Thus, if $t\mapsto \eta(x_0+tv)$ has positive second derivative at $t=0$, then $v^T D^2 \eta (x_0) v>0$. However, $v$ was arbitrary so this shows that $\eta$ is strictly convex at $x_0$.

We now show  $g(t)\coloneqq \eta(x_0+tv)$ has positive second derivative.

Write 
\begin{align}\label{finite_difference2nd}
g''(t) = \lim_{h\to0} \frac{g(t+h)-2g(t)+g(t-h)}{h^2}.
\end{align}
Then, note 
\begin{equation}
\begin{aligned}\label{convexity_breakdown}
& \frac{g(t+h)-2g(t)+g(t-h)}{h^2}\\
 &\hspace{.4in}=\frac{1}{h^2}\int\big[\eta_0(x_0+(t+h)v-y)-2\eta_0(x_0+tv-y)+\eta_0(x_0+(t-h)v-y)\big]m(y)\,dy.
 \end{aligned}
 \end{equation}
 
Then remark that the maximum which is taken in \eqref{max_for_convexity}, it is only over a finite index set. Thus, for each $x_0, t, v, y$, there exists $i$ such that $\eta_0(x_0+tv-y)=h_i+D_i\cdot(x_0+tv-y-x_i)+\epsilon_0\abs{x_0+tv-y-x_i}^2.$ For this $i$,
\begin{equation}
\begin{aligned}\label{convexity_math}
\eta_0(x_0+(t+h)v-y)-2\eta_0(x_0+tv-y)+\eta_0(x_0+(t-h)v-y)\hspace{1in}\\
\geq 
\Big(h_i+D_i\cdot(x_0+(t+h)v-y-x_i)+\epsilon_0\abs{x_0+(t+h)v-y-x_i}^2\Big) \\
-2\Big(h_i+D_i\cdot(x_0+tv-y-x_i)+\epsilon_0\abs{x_0+tv-y-x_i}^2\Big)\\
+\Big(h_i+D_i\cdot(x_0+(t-h)v-y-x_i)+\epsilon_0\abs{x_0+(t-h)v-y-x_i}^2\Big)\\
=2\epsilon_0 h^2 \norm{v}^2,
\end{aligned}
\end{equation}
where note that the last line is independent of $i$. From \eqref{finite_difference2nd}, \eqref{convexity_breakdown}, and \eqref{convexity_math}, we see that the map $\mathbb{R}\ni t\mapsto \eta(x_0+tv)$ has positive second derivative at $t=0$. 

This completes the proof of \Cref{eta_extension_prop}.

\section{Appendix: Proof of entropy conditions on shocks for systems from Baiti-Jenssen \cite{MR1849664}}\label{sec:laxliu_appendix}

The work \cite{MR1849664} does not explicitly discuss the Lax \emph{E}-condition \eqref{lax_cond} or Liu \emph{E}-condition \eqref{liu_cond}.

In this appendix, we show that the systems under consideration in \cite{MR1849664} have the property that each shock curve admits a family of shocks with both the Lax \emph{E}-condition \eqref{lax_cond} and Liu \emph{E}-condition \eqref{liu_cond}.

Our argument is based on a similar argument in \cite{K-K}.

The systems in \cite{MR1849664} have three characteristic fields. The middle field is simply a decoupled copy of a genuinely nonlinear scalar conservation law. This characteristic field will have a family of shocks verifying the Lax \emph{E}-condition and Liu \emph{E}-condition, simply due to the scalar theory. For Lax \emph{E}-condition for scalar, see \cite[p.~275]{dafermos_big_book}. For Liu \emph{E}-condition for scalar, see \cite[p.~279]{dafermos_big_book}.

The first and third fields are genuinely nonlinear with straight shock curves which coincide with the integral curves of the eigenvector fields. Remark that the first and third fields actually arise from a $2\times 2$ system which has coefficients depending on the solution to the decoupled scalar conservation law (the middle family of the $3\times3$ system). 

Assume we have parameterized the curves $S^1$ and $S^3$ with respect to arc length. Then, \cite[p.~456]{K-K} gives a differential equation which holds along each shock curve  in the Hugoniot locus $H(U)$ centered at point $U$ (stated for a $2\times 2$ system, but holding more generally):
\begin{align}\label{ode}
    \frac{\mathrm{d}}{\mathrm{d}s}\sigma^k_U(s) T=\alpha_1(\lambda_1-\sigma^k_U(s))r_1+\alpha_2(\lambda_2-\sigma^k_U(s))r_2+\alpha_3(\lambda_3-\sigma^k_U(s))r_3,
\end{align}
where at a point $S^k_U(s)$ in the Hugoniot locus, $T=S^k_U(s)-U$ and $t=\frac{\mathrm{d}}{\mathrm{d}s} T=\frac{\mathrm{d}}{\mathrm{d}s}S^k_U(s)$ is the unit tangent vector to the Hugoniot locus. Furthermore, in the context of \eqref{ode} we write 
\begin{align}\label{tangent}
    t=\alpha_1 r_1+\alpha_2 r_2+\alpha_3 r_3,
\end{align}
which is a representation of the tangent vector in terms of the local eigenvectors.

We will show now that the Lax \emph{E}-condition and Liu \emph{E}-condition hold along half of each of the two shocks curves $S^1$ and $S^3$ for this $3\times 3$ system (that is, for either $s<0$ or $s>0$). We follow an argument from \cite[p.~458-459]{K-K}.

Let us consider only the shock curve $S^1_U$ originating at $U$. The case for $S^3_U$ is nearly identical.

From \cite[Equations (8.2.1)-(8.2.2)]{dafermos_big_book}, we have that 
\begin{equation}
    \begin{aligned}\label{deriv_speed}
    \frac{\mathrm{d}}{\mathrm{d}s}\sigma^1_U(0)&=\frac{1}{2}\nabla\lambda_1(U) \cdot r_1(U),\\
    \sigma^1_U(0)&=\lambda_1(U).
\end{aligned}
\end{equation}

From \eqref{deriv_speed} and genuine nonlinearity, both the Liu \emph{E}-condition and Lax \emph{E}-condition will hold in a small (one-sided) neighborhood $N$ of $s=0$ on $S^1_U$.

Moreover, from \eqref{deriv_speed} and genuine nonlinearity,
\begin{align}\label{strong_liu}
\frac{\mathrm{d}}{\mathrm{d}s}\sigma^1_U\neq 0
\end{align}
will hold in the (one-sided) neighborhood $N$ of $s=0$ on $S^1_U$ (possibly shrinking the size of $N$ if necessary). Note that \eqref{strong_liu} is stronger than the Liu \emph{E}-condition.
  
By again possibly shrinking $N$ if necessary,
\begin{align}\label{strong_lax}
    \lambda_1(S^1_{U}(s))< \sigma^1_{U}(s) < \lambda_1(U),
\end{align}
will also will hold along $N$. Remark that \eqref{strong_lax} is stronger than the Lax \emph{E}-condition.

Without loss of generality, we can assume that $S^1_U$ with $s\geq0$ contains the half-neighborhood $N$ (and thus, $\frac{\mathrm{d}}{\mathrm{d}s}\sigma^1_U(0)< 0$).

We now show that the neighborhood $N$ where \eqref{strong_liu} and \eqref{strong_lax} holds actually contains all of $S^1_U$, for $s\geq0$. If it does not, then there will be a first point $U_1\neq U$ on $S^1_U$ at which \eqref{strong_liu} or one of the inequalities \eqref{strong_lax} is violated.

The second half of \eqref{strong_lax} cannot be violated at $U_1$ because \eqref{strong_liu} holds in the rest of $N$.

Remark that along $S^1_U$, we have $\alpha_1=1$ and $\alpha_2=\alpha_3=0$ in \eqref{tangent}, because the shock curve coincides with the integral curve of the eigenvector field.

Then, from \eqref{ode} we must have that at this point of violation $U_1$, $\frac{\mathrm{d}}{\mathrm{d}s}\sigma^1_U=0$ and $\sigma^1_U=\lambda_1$. 

However, from genuine nonlinearity we have $\frac{\mathrm{d}}{\mathrm{d}s} (\sigma^1_U(s)-\lambda_1(S^1_U(s)))>0$ at $U_1$ due to  $\frac{\mathrm{d}}{\mathrm{d}s}\sigma^1_U(0)< 0$ (recalling also that the shock curve coincides with the integral curve of the corresponding eigenvector field). But then, due to $\sigma^1_U-\lambda_1=0$ at $U_1$, we must have $\sigma^1_U(s)-\lambda_1(S^1_U(s))<0$ for points of $N$ near $U_1$. This contradicts \eqref{strong_lax}. Thus, $N$ must equal all of $S^1_U$ (for $s\geq0$).

\section{Appendix: Large $T_5$ example}\label{sec:T5_example}

We present an example of a large $T_5$, from the context of \Cref{matlab_lemma}.

\begin{align}
    X_1^0 = 
 \begin{pmatrix}
  \frac{2316901181183546099017091770605}{162259276829213363391578010288128} & \frac{-5739565125354385225872252139467}{10141204801825835211973625643008} \\[0.9em]
  \frac{5739565125354385225872252139467}{10141204801825835211973625643008} & \frac{-56873565598434451746179265262997}{2535301200456458802993406410752}
 \end{pmatrix}
\end{align}

\begin{align}
    X_2^0 = 
 \begin{pmatrix}
  \frac{-330621418565185387036477002414115}{10384593717069655257060992658440192} & \frac{815211547324287408551802829225453}{649037107316853453566312041152512} \\[0.9em]
  \frac{-815211547324287408551802829225453}{649037107316853453566312041152512} & \frac{251257742411123530141636860664321}{5070602400912917605986812821504}
 \end{pmatrix}
\end{align}

\begin{align}
    X_3^0 = 
 \begin{pmatrix}
  \frac{1231069874758438218672166401101419}{20769187434139310514121985316880384} & \frac{-733625671232943434981364913268293}{324518553658426726783156020576256} \\[0.9em]
  \frac{733625671232943434981364913268293}{324518553658426726783156020576256} & \frac{-3497495142386849315118349227834509}{40564819207303340847894502572032}
 \end{pmatrix}
\end{align}

\begin{align}
    X_4^0 = 
 \begin{pmatrix}
  \frac{-587406988058843286220046809310939}{20769187434139310514121985316880384} & \frac{88287046124795489454709265219111}{81129638414606681695789005144064} \\[0.9em]
  \frac{-88287046124795489454709265219111}{81129638414606681695789005144064} & \frac{1698486469749761796168679052449441}{40564819207303340847894502572032}
 \end{pmatrix}
\end{align}

\begin{align}
    X_5^0 = 
 \begin{pmatrix}
  \frac{-907999771200015425284613822879209}{41538374868278621028243970633760768} & \frac{1065154343959802482930848455328023}{1298074214633706907132624082305024} \\[0.9em]
  \frac{-1065154343959802482930848455328023}{1298074214633706907132624082305024} & \frac{2499017758464245906295246824575579}{81129638414606681695789005144064}
 \end{pmatrix}
\end{align}

\bibliographystyle{plain}
\bibliography{references}

\begin{thebibliography}{10}

\bibitem{fmincon}
{F}ind minimum of constrained nonlinear multivariable function -- {MATLAB}
  fmincon.
\newblock \url{https://www.mathworks.com/help/optim/ug/fmincon.html}.
\newblock Accessed: March 3, 2024.

\bibitem{MATLAB_double}
{F}loating-{P}oint {N}umbers - {MATLAB} \& {S}imulink.
\newblock
  \url{https://www.mathworks.com/help/matlab/matlab_prog/floating-point-numbers.html}.
\newblock Accessed: March 10, 2024.

\bibitem{MATLAB_solve}
{S}olve optimization problem or equation problem -- {MATLAB} solve.
\newblock
  \url{https://www.mathworks.com/help/optim/ug/optim.problemdef.optimizationproblem.solve.html}.
\newblock Accessed: March 3, 2024.

\bibitem{MATLAB_symbolic}
{S}ymbolic {M}ath {T}oolbox {D}ocumentation.
\newblock \url{https://www.mathworks.com/help/symbolic/}.
\newblock Accessed: March 3, 2024.

\bibitem{MR1849664}
Paolo Baiti and Helge~Kristian Jenssen.
\newblock Blowup in {$L^\infty$} for a class of genuinely nonlinear hyperbolic
  systems of conservation laws.
\newblock {\em Discrete Contin. Dynam. Systems}, 7(4):837--853, 2001.

\bibitem{MR2150387}
Stefano Bianchini and Alberto Bressan.
\newblock Vanishing viscosity solutions of nonlinear hyperbolic systems.
\newblock {\em Ann. of Math. (2)}, 161(1):223--342, 2005.

\bibitem{MR1816648}
Alberto Bressan.
\newblock {\em Hyperbolic systems of conservation laws}, volume~20 of {\em
  Oxford Lecture Series in Mathematics and its Applications}.
\newblock Oxford University Press, Oxford, 2000.
\newblock The one-dimensional Cauchy problem.

\bibitem{MR3914001}
Alberto Bressan, Geng Chen, and Qingtian Zhang.
\newblock On finite time {BV} blow-up for the p-system.
\newblock {\em Comm. Partial Differential Equations}, 43(8):1242--1280, 2018.

\bibitem{MR4661213}
Alberto Bressan and Camillo De~Lellis.
\newblock A remark on the uniqueness of solutions to hyperbolic conservation
  laws.
\newblock {\em Arch. Ration. Mech. Anal.}, 247(6):Paper No. 106, 12, 2023.

\bibitem{MR1723032}
Alberto Bressan, Tai-Ping Liu, and Tong Yang.
\newblock {$L^1$} stability estimates for {$n\times n$} conservation laws.
\newblock {\em Arch. Ration. Mech. Anal.}, 149(1):1--22, 1999.

\bibitem{MR4487515}
Geng Chen, Sam~G. Krupa, and Alexis~F. Vasseur.
\newblock Uniqueness and weak-{BV} stability for {$2\times 2$} conservation
  laws.
\newblock {\em Arch. Ration. Mech. Anal.}, 246(1):299--332, 2022.

\bibitem{multi_d_illposed}
Elisabetta Chiodaroli, Camillo De~Lellis, and Ond\v{r}ej Kreml.
\newblock Global ill-posedness of the isentropic system of gas dynamics.
\newblock {\em Comm. Pure Appl. Math.}, 68(7):1157--1190, 2015.

\bibitem{dafermos_big_book}
Constantine~M. Dafermos.
\newblock {\em Hyperbolic conservation laws in continuum physics}, volume 325
  of {\em Grundlehren der Mathematischen Wissenschaften [Fundamental Principles
  of Mathematical Sciences]}.
\newblock Springer-Verlag, Berlin, fourth edition, 2016.

\bibitem{MR4627977}
Constantine~M. Dafermos.
\newblock Hyperbolic conservation laws: past, present, future.
\newblock In {\em Mathematics going forward---collected mathematical
  brushstrokes}, volume 2313 of {\em Lecture Notes in Math.}, pages 479--486.
  Springer, Cham, [2023] \copyright 2023.

\bibitem{MR2600877}
Camillo De~Lellis and L\'{a}szl\'{o} Sz\'{e}kelyhidi, Jr.
\newblock The {E}uler equations as a differential inclusion.
\newblock {\em Ann. of Math. (2)}, 170(3):1417--1436, 2009.

\bibitem{MR2564474}
Camillo De~Lellis and L\'{a}szl\'{o} Sz\'{e}kelyhidi, Jr.
\newblock On admissibility criteria for weak solutions of the {E}uler
  equations.
\newblock {\em Arch. Ration. Mech. Anal.}, 195(1):225--260, 2010.

\bibitem{MR3740399}
Clemens F\"{o}rster and L\'{a}szl\'{o} Sz\'{e}kelyhidi, Jr.
\newblock {$T_5$}-configurations and non-rigid sets of matrices.
\newblock {\em Calc. Var. Partial Differential Equations}, 57(1):Paper No. 19,
  12, 2018.

\bibitem{MR0194770}
James~G. Glimm.
\newblock Solutions in the large for nonlinear hyperbolic systems of equations.
\newblock {\em Comm. Pure Appl. Math.}, 18:697--715, 1965.

\bibitem{MR864505}
Mikhael Gromov.
\newblock {\em Partial differential relations}, volume~9 of {\em Ergebnisse der
  Mathematik und ihrer Grenzgebiete (3) [Results in Mathematics and Related
  Areas (3)]}.
\newblock Springer-Verlag, Berlin, 1986.

\bibitem{MR1912206}
Helge Holden and Nils~Henrik Risebro.
\newblock {\em Front tracking for hyperbolic conservation laws}, volume 152 of
  {\em Applied Mathematical Sciences}.
\newblock Springer-Verlag, New York, 2002.

\bibitem{MR1752421}
Helge~Kristian Jenssen.
\newblock Blowup for systems of conservation laws.
\newblock {\em SIAM J. Math. Anal.}, 31(4):894--908, 2000.

\bibitem{https://doi.org/10.48550/arxiv.2208.10979}
Carl Johan~Peter {Johansson} and Riccardo {Tione}.
\newblock {$T_5$} configurations and hyperbolic systems.
\newblock {\em Communications in Contemporary Mathematics}, 0(0):2250081, 2023.

\bibitem{MR1272766}
Jean-Luc Joly, Guy M\'{e}tivier, and Jeffrey Rauch.
\newblock A nonlinear instability for {$3\times 3$} systems of conservation
  laws.
\newblock {\em Comm. Math. Phys.}, 162(1):47--59, 1994.

\bibitem{K-K}
Barbara~L. Keyfitz and Herbert~C. Kranzer.
\newblock Existence and uniqueness of entropy solutions to the {R}iemann
  problem for hyperbolic systems of two nonlinear conservation laws.
\newblock {\em J. Differential Equations}, 27(3):444--476, 1978.

\bibitem{hab_thesis}
Bernd Kirchheim.
\newblock {\em Rigidity and Geometry of microstructures}.
\newblock 2003.
\newblock Habilitation thesis -- University of Leipzig.

\bibitem{MR2008346}
Bernd Kirchheim, Stefan M\"{u}ller, and Vladim\'{\i}r \v{S}ver\'{a}k.
\newblock Studying nonlinear {PDE} by geometry in matrix space.
\newblock In {\em Geometric analysis and nonlinear partial differential
  equations}, pages 347--395. Springer, Berlin, 2003.

\bibitem{2022arXiv221114239K}
Sam~G. {Krupa} and L\'{a}szl\'{o} Sz\'{e}kelyhidi, Jr.
\newblock {Nonexistence of $T_4$ configurations for hyperbolic systems and the
  Liu entropy condition}.
\newblock {\em arXiv e-prints}, page arXiv:2211.14239, November 2022.

\bibitem{MR4144350}
Andrew Lorent and Guanying Peng.
\newblock On the rank-1 convex hull of a set arising from a hyperbolic system
  of {L}agrangian elasticity.
\newblock {\em Calc. Var. Partial Differential Equations}, 59(5):Paper No. 156,
  36, 2020.

\bibitem{MR4385531}
Simon Markfelder.
\newblock {\em Convex integration applied to the multi-dimensional compressible
  {E}uler equations}, volume 2294 of {\em Lecture Notes in Mathematics}.
\newblock Springer, Cham, [2021] \copyright 2021.

\bibitem{MR1728376}
S.~M\"{u}ller and V.~\v{S}ver\'{a}k.
\newblock Convex integration with constraints and applications to phase
  transitions and partial differential equations.
\newblock {\em J. Eur. Math. Soc. (JEMS)}, 1(4):393--422, 1999.

\bibitem{MR1983780}
Stefan M\"{u}ller and Vladim\'{\i}r \v{S}ver\'{a}k.
\newblock Convex integration for {L}ipschitz mappings and counterexamples to
  regularity.
\newblock {\em Ann. of Math. (2)}, 157(3):715--742, 2003.

\bibitem{MR2624766}
Vladimir Scheffer.
\newblock {\em Regularity and irregularity of solutions to nonlinear second
  order elliptic systems and inequalities}.
\newblock ProQuest LLC, Ann Arbor, MI, 1974.
\newblock Thesis (Ph.D.)--Princeton University.

\bibitem{MR1301779}
Joel~Alan Smoller.
\newblock {\em Shock waves and reaction-diffusion equations}, volume 258 of
  {\em Grundlehren der mathematischen Wissenschaften [Fundamental Principles of
  Mathematical Sciences]}.
\newblock Springer-Verlag, New York, second edition, 1994.

\bibitem{MR2048569}
L\'{a}szl\'{o} Sz\'{e}kelyhidi, Jr.
\newblock The regularity of critical points of polyconvex functionals.
\newblock {\em Arch. Ration. Mech. Anal.}, 172(1):133--152, 2004.

\bibitem{MR2118899}
L\'{a}szl\'{o} Sz\'{e}kelyhidi, Jr.
\newblock Rank-one convex hulls in {$\Bbb R^{2\times 2}$}.
\newblock {\em Calc. Var. Partial Differential Equations}, 22(3):253--281,
  2005.

\bibitem{MR1320538}
Luc Tartar.
\newblock Some remarks on separately convex functions.
\newblock In {\em Microstructure and phase transition}, volume~54 of {\em IMA
  Vol. Math. Appl.}, pages 191--204. Springer, New York, 1993.

\bibitem{MR2508169}
Alexis~F. Vasseur.
\newblock Recent results on hydrodynamic limits.
\newblock In {\em Handbook of differential equations: evolutionary equations.
  {V}ol. {IV}}, Handb. Differ. Equ., pages 323--376. Elsevier/North-Holland,
  Amsterdam, 2008.

\bibitem{MR1677943}
Robin Young.
\newblock Exact solutions to degenerate conservation laws.
\newblock {\em SIAM J. Math. Anal.}, 30(3):537--558, 1999.

\end{thebibliography}
\end{document}